\documentclass[12pt,a4paper]{article}
\usepackage[english]{babel}
\usepackage[OT1]{fontenc}

\usepackage{amsmath, amssymb}
\usepackage[bf]{caption2}
\usepackage{amsfonts}
\usepackage{amssymb}
\usepackage{graphicx}
\usepackage{graphics}
\usepackage[usenames]{color}
\usepackage{colortbl}
\usepackage{dblfloatfix}
\usepackage{fixltx2e}
\usepackage{mdwtab}

\newcounter{theorem}

\newcommand{\la}{\lambda}
\newcommand{\vp}{{\mathbf p}}
\newcommand{\vz}{{\mathbf z}}

\newcommand{\vf}{{\mathbf f}}
\newcommand{\vx}{{\mathbf x}}

\newcommand{\be}{\begin{equation}}
\newcommand{\ee}{\end{equation}}
\newcommand{\ber}{\begin{eqnarray}}
\newcommand{\eer}{\end{eqnarray}}

\begin{document}

\title{On limiting characteristics  for  a non-stationary two-processor heterogeneous system}

\author{A. Zeifman\footnote{Vologda State University; Institute of Informatics
Problems, Federal Research Center "Informatics and Control" of RAS; Vologda Research Center of RAS; $a\_$zeifman@mail.ru}, Y. Satin\footnote{Vologda State University},  K. Kiseleva\footnote{Vologda State University}, T. Panfilova\footnote{Vologda State University}, V. Korolev\footnote{Faculty of
Computational Mathematics and Cybernetics, Lomonosov Moscow State
University; Institute of Informatics Problems, Federal Research
Center "Informatics and Control" of RAS}}

\date{}

\maketitle

\section{Introduction}

In this paper we study a non-stationary Markovian queueing model of
a two-processor heterogeneous system with time-varying arrival and
service rates which was firstly investigated in \cite{dh2000}, see
also time-dependent analysis of this model in the recent paper
\cite{aa2018}. In general, non-stationary queueing models have been
actively studied during some decades, see, for instance
\cite{Dong2015,Giorno2014,Granovsky2004,Schwarz2016,Zeifman2006,Zeifman2014q}
and the references therein.

In the paper  \cite{aa2018} the authors deal with the so-called
``time-dependent analysis'', in other words, they try to find the
state probabilities on a finite interval under some initial
conditions (as a rule, initially, the number of customers in the
queue is zero), see for instance \cite{DiCrescenzo2016}. Another
approach is connected with the determination of the limiting mode,
see \cite{Chakravarthy2017}.

Essentially more information about queue-length process can be
obtained using ergodicity and the corresponding estimates of the
rate of convergence. A general approach to obtaining sharp bounds on
the rate of convergence via the notion of the logarithmic norm of an
operator function wsa recently discussed in details in our papers
\cite{Zeifman2018jamcs,Zeifman2018lncs,Zeifman2018spl}.  The first
studies in this direction were published since 1980-s for
birth-death models, see \cite{Zeifman1989ait,Zeifman1995s}. In
\cite{Zeifman2018jamcs} we proved that there are four classes of
Markovian queueing models for which the reduced forward Kolmogorov
system can be transformed to the system with essentially nonnegative
matrix. Although the model under consideration does not belong to
one of these classes, we can apply the same approach and obtain some
useful bounds on the rate of convergence for it. Moreover, we can
compute the limiting characteristics of the model using bounds on
the rate of convergence and truncations technique introduced in
\cite{Zeifman2014i,Zeifman2016t}.

Note an interesting fact: exact estimates of the rate of convergence
yield exact estimates of stability (perturbation bounds), see
\cite{Kartashov1985,Liu2012,Mitrophanov2003,Mitrophanov2004,Zeifman1985,Zeifman2014s}
and references therein.

An important feature of multiprocessor queueing systems is the
presence of risks related to the overload of the system. In the
present paper it is demonstrated that under natural conditions on
the arrival/service rates these risks vanish and the system rather
easily approaches the ergodic mode.

\section{Description of the model}

Here we consider a multiprocessor system consisting of two types of
processors, which for convenience will be referred to as the
``main'' and ``backup'' processors \cite{aa2018}. Each job requires
exactly one processor for its execution. When both processors are
idle, the main processor is scheduled for service before the backup
processor. A computer system consists of two processors, a main
processor, and a backup processor. A description of the model is as
follows:

\smallskip
\noindent (i) jobs arrive at the system according to the Poisson
process with an arrival rate $\lambda(t)$. Service is exponentially
distributed, and two servers provide heterogeneous service rates
$\mu_1(t),~\mu_2(t)$ such that $\mu_2(t)\le\mu_1(t)$.

\smallskip
\noindent (ii) each job needs only one server to be served and the
jobs select the servers on the basis of fastest server first (FSF).

\begin{figure}[htbp]
\begin{center}
  \includegraphics[scale=0.31]{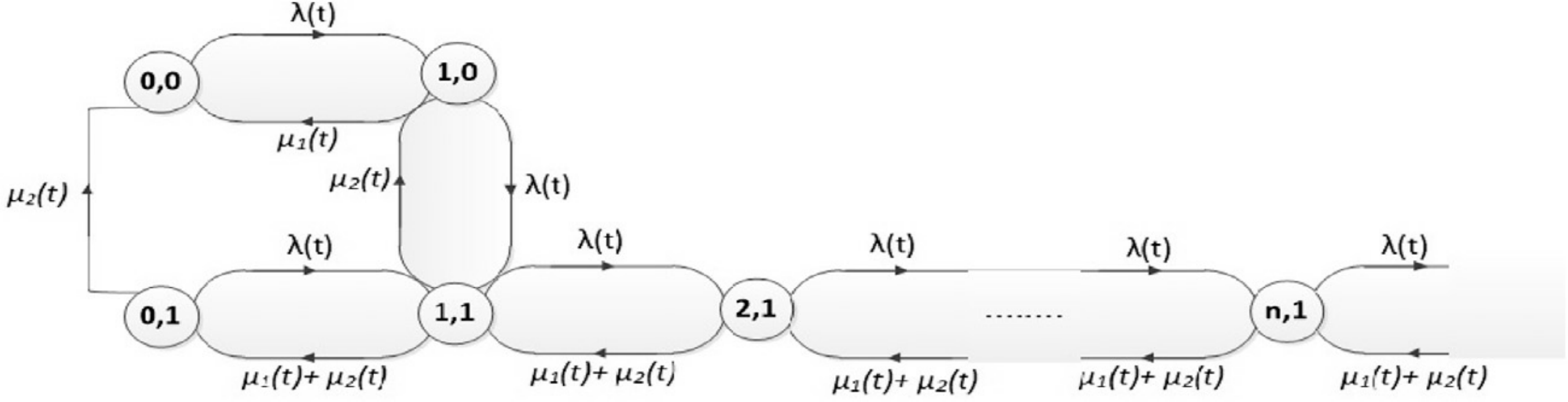}
\end{center}
\caption{Transitions for a two-processor heterogeneous model}
\label{Descr}
\end{figure}

The probabilistic dynamics of the process is represented by the
forward Kolmogorov system of differential equations:
\begin{equation}
\frac{d\vp}{dt}=A(t)\vp,
\label{mains}
\end{equation}
where $\vp=(p_{00}, p_{10}, p_{01}, p_{11}, p_{12}, \ldots, p_{1n},\ldots)^T,$

\begin{equation}
A(t)=\left(
\begin{array}{cccccccc}
-\la & \mu_1  & \mu_2  & 0 & 0 &\cdots \\
\la  & -(\la+\mu_1) & 0& \mu_2 & 0 & \cdots \\
0  & 0 & -(\la+\mu_2)& \mu_1 & 0  & \cdots \\
0  & \la & \la&  -(\la+\mu)& \mu   & \cdots \\
0  & 0 & 0 & \la&  -(\la+\mu)  & \cdots \\
0  & 0 & 0 & 0 &  \la &  \cdots \\
\cdots & \cdots & \cdots & \cdots & \cdots \\
\end{array}
\right),
\end{equation}
where $\mu(t)=\mu_1(t)+\mu_2(t),~A(t)=Q^T(t)$, and $Q(t)$ - the intensity matrix.

\section{Bounds on the rate of convergence}
Since $p_{00}(t) = 1 - p_{01}(t)-\sum_{j=0}^{\infty} p_{1j}(t)$ due
to the normalization condition, the system (\ref{mains}) can be
rewritten as
\begin{equation}
\frac{d \vz}{dt}= B(t)\vz+\vf(t), \label{sys2}
\end{equation}
\noindent where
$$
\vf(t)=\left( \la,  0, \ldots, 0, \ldots\right)^{T}, \
\vz(t)=\left(p_{10}, p_{01}, p_{11}, p_{12}, \ldots, p_{1n},\ldots \right)^{T},$$ and
$B(t)=(b_{ij}(t))_{i,j=1}^{\infty}=$

\begin{equation}
=\left(
\begin{array}{cccccccc}
-(2\la+\mu_1) & -\la  & \mu_2-\la  & -\la & -\la &\cdots \\
0  & -(\la+\mu_2) & \mu_1& 0 & 0 & \cdots \\
\la  & \la & -(\la+\mu)& \mu & 0  & \cdots \\
0  & 0 & \la & -(\la+\mu) &  \mu  & \cdots \\
\cdots & \cdots & \cdots & \cdots & \cdots \\
\end{array}
\right).
\end{equation}

Denote by $T$ the upper triangular matrix
\begin{equation}
T=\left(
\begin{array}{ccccccc}
1   & 1 & 1 & \cdots & 1 \\
0   & 1  & 1  &   \cdots & 1 \\
0   & 0  & 1  &   \cdots & 1 \\
0   & 0  & 0  &   \cdots & 1\\
\vdots & \vdots & \vdots & \ddots \\
\end{array}
\right). \label{vspmatr}
\end{equation}

Consider the matrix $T B(t)T^{-1}=$
\begin{equation}
=\left(
\begin{array}{cccccccc}
-(\la+\mu_1) & \mu_1-\mu_2  & \mu_2  & 0 & 0 &\cdots \\
\la  & -(\la+\mu_2) & 0 & \mu_2 & 0 & \cdots \\
\la  & 0 & -(\la+\mu)& \mu & 0  & \cdots \\
0  & 0 & \la & -(\la+\mu) &  \mu  & \cdots \\
\cdots & \cdots & \cdots & \cdots & \cdots \\
\end{array}
\right).
\end{equation}

Let  $\{d_i\}$, $i \ge 1$, be a sequence of positive numbers such
that
$$d_1=1,~d_2=\epsilon,~d_3=1,~d_4=\delta_1>1, ~~\frac{d_5}{d_4}=\frac{d_6}{d_5}=\ldots=\delta >1.$$
Let $D = diag\left(d_1, d_2, \dots \right)$ be the corresponding
diagonal matrix and $l_{1\textsf{D}}$ be a space of vectors
$l_{1\textsf{D}}=\left\{{\bf x} =(x_1,x_2,\ldots)/\|{\bf
x}\|_{1\textsf{D}}=\|\textsf{D} {\bf x}\|_1 <\infty\right\},$ where
$\textsf{D}=DT.$

Consider the matrix $DT
B(t)T^{-1}D^{-1}=\textsf{D}B(t)\textsf{D}^{-1}=$ {\small
\begin{equation} =\left(
\begin{array}{ccccccccc}
-(\la+\mu_1) & \frac{d_1}{d_2}(\mu_1-\mu_2)  & \frac{d_1}{d_3}\mu_2  & 0 & 0 & 0 &\cdots \\
\\
\frac{d_2}{d_1}\la  & -(\la+\mu_2) & 0 & \frac{d_2}{d_4}\mu_2 & 0 & 0 &\cdots \\
\\
\frac{d_3}{d_1}\la  & 0 & -(\la+\mu)& \frac{d_3}{d_4}\mu & 0  & 0 & \cdots \\
\\
0  & 0 & \frac{d_4}{d_3}\la & -(\la+\mu) &  \frac{d_4}{d_5}\mu & 0 & \cdots \\
\\
0  & 0 & 0 & \frac{d_5}{d_4}\la & -(\la+\mu) &  \frac{d_5}{d_6}\mu  & \cdots \\
\\
0  & 0 & 0 & 0 & \frac{d_6}{d_5}\la &  -(\la+\mu)  & \cdots \\
\\
0  & 0 & 0 & 0 & 0 & \frac{d_7}{d_6}\la  & \cdots \\
\cdots & \cdots & \cdots & \cdots & \cdots & \cdots \\
\end{array}
\right).
\label{finmatrix}
\end{equation}}

The approach used in this paper is based on the notion of the
logarithmic norm of a linear operator function and the corresponding
bounds of the Cauchy operator, see the detailed discussion, for
instance, in \cite{Doorn2009}. Namely, if $B\left( t\right)$, $t\ge
0$, is a one-parameter family of bounded linear operators on a
Banach space ${\cal B}$, then
\begin{eqnarray}
\gamma \left( B\left( t\right) \right)_{\cal B}
=\lim\limits_{h\rightarrow
+0}\frac{%
\left\| I+hB\left( t\right) \right\| -1}h \label{lognorm}
\end{eqnarray}
is called the logarithmic norm of the operator $B\left( t\right)$.

If ${\cal B}=l_1$, then  the operator $B\left(t\right)$ is given by
the matrix $B\left(t\right)=\left( b_{ij}\left(t\right)\right)
_{i,j=0}^\infty$, $t\ge 0$, and the logarithmic norm of
$B\left(t\right)$ can be found explicitly:
\begin{equation}
\gamma(B(t))_{1\textsf{D}}=\gamma(\textsf{D}B(t)\textsf{D}^{-1})_1=\sup\limits_j\bigg(
b_{jj}\left( t\right) +\sum\limits_{i\neq j}\left| b_{ij}\left(
t\right) \right| \bigg) ,\quad t\ge 0.
\end{equation}

Hence the following bound on the rate of convergence holds:
$$ \|\vx
(t)\| \le e^{\int_0^t \gamma \left( B\left( \tau \right) \right)\,
d\tau}\|\vx(0)\|, $$ where $\vx=\vz^*-\vz^{**}$ and $\vx$ is the solution of the
differential equation
$$
\frac{d \vx}{dt}= B(t)\vx,
$$
which we obtain instead of the system (\ref{sys2}).

Let $\alpha_i(t)$ be negative sums of the elements of corresponding
columns for the matrix (\ref{finmatrix}), such as:
$$\alpha_1=(\la+\mu_1)-\frac{d_2}{d_1}\la-\frac{d_3}{d_1}\la, $$
$$\alpha_2=(\la+\mu_2)-\frac{d_1}{d_2}(\mu_1-\mu_2),$$
$$\alpha_3=(\la+\mu)-\frac{d_1}{d_3}\mu_2-\frac{d_4}{d_3}\la,$$
$$ \alpha_4=(\la+\mu)-\frac{d_2}{d_4}\mu_2-\frac{d_3}{d_4}\mu-\frac{d_5}{d_4}\la, $$
$$\alpha_5=(\la+\mu)-\frac{d_4}{d_5}\mu-\frac{d_6}{d_5}\la, ~~~ \alpha_6=(\la+\mu)-\frac{d_5}{d_6}\mu-\frac{d_7}{d_6}\la,~~~ \ldots,$$
where $\alpha_5=\alpha_6=\ldots,$ since
$\frac{d_5}{d_4}=\frac{d_6}{d_5}=\frac{d_7}{d_6}=\ldots=\delta$.

Then we obtain the logarithmic norm:
\begin{equation}
\gamma(B(t))_{1\textsf{D}}=\gamma(\textsf{D}B(t)\textsf{D}^{-1})_1=-\inf_{i \ge 1}(\alpha_i(t)) = -\min_{i \le 5}(\alpha_i(t)).
\end{equation}

\section{The case $\mu_1 =\mu_2$.}

First, let $\la,~\mu_1 = \mu_2$ be constant,
$0<\la<\mu=\mu_1+\mu_2.$ Then the exact value of the decay parameter
(or the spectral gap) for a simple birth-death process with
intensities $\la$ and $\mu$ is well-known, namely, it equals
$\beta^*=(\sqrt{\la}-\sqrt{\mu})^2,$ see, e. g., \cite{Doorn1985},
and the corresponding $\delta=\sqrt{\frac{\mu}{\la}}~$. Hence, we
consider the same $\delta$ and put $d_2=\epsilon<<1,$
$\delta_1=\delta$.

Then
$$\alpha_1=\frac{\mu}{2}-\epsilon\la,$$
$$\alpha_2=\la+\frac{\mu}{2},$$
$$\alpha_3=\frac{\mu}{2}+\la-\sqrt{\la\mu},$$
$$\alpha_4=(\sqrt{\la}-\sqrt{\mu})^2-\frac{\epsilon}{2}\sqrt{\la\mu}, $$
$$\alpha_k=(\sqrt{\la}-\sqrt{\mu})^2,~~~k\ge5,$$

Put $\beta_*=\min(\alpha_1,~\alpha_3,~\alpha_4).$ Then we have
\begin{equation}
\gamma(B(t))_{1\textsf{D}}=-\inf_i \alpha_i(t) =-\beta_*(t).
\end{equation}
Hence the following bound holds:
\begin{equation}
\|\vx(t)\|_{1\textsf{D}} \le e^{-\beta_*t}\|\vx(0)\|_{1\textsf{D}}.
\end{equation}

Let now the intensities $\la(t),~\mu_1(t) = \mu_2(t)$ be 1-periodic.
Put
$$\mu_*=\int_0^1\mu(t)dt,~~\la_*=\int_0^1\la(t)dt.$$
Then the best possible bound for a pure birth-death process is
attained, if we take $\delta=\sqrt{\frac{\mu_*}{\la_*}}.$

Then for these $\delta$, $d_2=\epsilon<<1,$ $\delta_1=\delta$ we
have
\begin{equation}
\alpha_1(t)=\frac{\mu(t)}{2}-\epsilon\la(t), \quad
\alpha_2(t)=\la(t)+\frac{\mu(t)}{2}, \label{alpha01}
\end{equation}
\begin{equation}\alpha_3(t)=\frac{\mu(t)}{2}+\la(t)-\sqrt{\la(t)\mu(t)},
\label{alpha03}
\end{equation}
\begin{equation}\alpha_4(t)=(\sqrt{\la(t)}-\sqrt{\mu(t)})^2-\frac{\epsilon}{2}\sqrt{\la(t)\mu(t)},
\label{alpha04}
\end{equation}
\begin{equation}
\alpha_k(t)=(\sqrt{\la(t)}-\sqrt{\mu(t)})^2,~~~k\ge 5,
\label{alpha05}
\end{equation}

Put $\beta_*(t)=\min(\alpha_1(t),~\alpha_3(t),~\alpha_4(t)).$ We
have
\begin{equation}
\gamma(B(t))_{1\textsf{D}}=-\inf \alpha_i(t) = -\beta_*(t).
\label{beta01}
\end{equation}
Hence the following bound holds:
\begin{equation}
\|\vx(t)\|_{1\textsf{D}} \le e^{-\int_0^t\beta_*(\tau)\,dt}\|\vx(0)\|_{1\textsf{D}}.
\end{equation}

\smallskip

{\bf Remark 1.} It should be noted that in \cite{aa2018} there are
some misprints in the plots, namely, the intensities must have a
multiplier $\pi$, say $1+\sin2 \pi t$. Moreover, on Fig 3 of that
paper the sum of all probabilities evidently is greater than 1.

\smallskip

{\bf Remark 2.} It can be seen that actually the periodic terms in
the intensities do not affect the rate of convergence, see the plots
related to the examples. Hence it is essentially easier to find the
parameter $\beta_{0*}$ for the corresponding homogeneous model.
Namely, if we put
\begin{equation}
\alpha_{10}=\frac{\mu_*}{2}-\epsilon\la_*, \quad \alpha_{20}=\la_*
+\frac{\mu_*}{2},
\alpha_{30}=\frac{\mu_*}{2}+\la_*-\sqrt{\la_*\mu_*},
\label{alpha*01}
\end{equation}
\begin{equation}\alpha_{40}=(\sqrt{\la_*}-\sqrt{\mu_*})^2-\frac{\epsilon}{2}\sqrt{\la_*\mu_*},
 \alpha_{k0}=(\sqrt{\la_*}-\sqrt{\mu_*})^2,~~~k\ge 5,
\label{alpha*05}
\end{equation}
instead of (\ref{alpha01})-(\ref{alpha05}), then we obtain
$\beta_{*0}=\min\left(\alpha_{10},~\alpha_{30},~\alpha_{40}\right),$
and the bound on the rate of convergence in the form
\begin{equation}
\|\vx(t)\|_{1\textsf{D}} \le \textsf{N}e^{-\beta_{*0}
t}\|\vx(0)\|_{1\textsf{D}}, \label{beta*01}
\end{equation}
for some positive $\textsf{N}$.

\smallskip

{\bf Remark 3.} In the following examples we consider the behavior
of the 'first' state probabilities $P_{00}(t)$, $P_{01}(t)$,
$P_{10}(t)$, $P_{11}(t)$, $P_{21}(t)$, $P_{31}(t)$ and the
mathematical expectation (the mean) of the queue-length process
$E(t) = p01+p10+2*p11+3*p21+4*p31+ \dots$. One can see that for {\it
all} examples the rates of convergence for original model with
$1-$periodic intensities and for the corresponding homogeneous model
are the same, as it was noted in Remark 2.

\section{Examples}

{\bf Example 1.} Let $\mu_1=\mu_2=2,~~\la=1+\sin2\pi t, ~~\la_*=1$
(Example 2 from \cite{aa2018}). Put
$\delta=\sqrt{\frac{\mu}{\la_*}}=2$. Then we obtain
$$\alpha_1(t)=2-\epsilon(1+\sin2\pi t) \ge 2 - 2\epsilon,$$
$$\alpha_3=3+\sin2\pi t-2\sqrt{1+\sin2\pi t}= 1+\left(1-\sqrt{1+\sin2\pi t} \right)^2 \ge 1,$$
$$\alpha_4=\left(2- \sqrt{1+\sin2\pi t}\right)^2-\epsilon\sqrt{1+\sin2\pi t} \ge \left(2- \sqrt{2}\right)^2-\epsilon\sqrt{2} \ge 0.3,$$
for sufficiently small $\epsilon$. Therefore, we have $\beta_* =
0.3$. Thus, we can obtain the following {\bf bound}
\begin{equation*}
\|\vp^*(t)-\vp^{**}(t)\|_1\le2\|\vz^*(t)-\vz^{**}(t)\|_1\le4\|\vz^*(t)-\vz^{**}(t)\|_{1\textsf{D}}\le 4e^{-0.3 t}\|\vz^*(0)-\vz^{**}(0)\|_{1\textsf{D}}.
\end{equation*}

On the other hand, we can obtain a simpler bound by applying Remark
2. Namely, $ \alpha_{10}= 2-\epsilon$, $ \alpha_{30}=1$, $
\alpha_{40}=1-\epsilon$, $\beta_{*0}= 1-\epsilon$, and
\begin{equation*}
\|\vp^*(t)-\vp^{**}(t)\|_1  \le
4\textsf{N}e^{-\left(1-\epsilon\right)
t}\|\vz^*(0)-\vz^{**}(0)\|_{1\textsf{D}}.
\end{equation*}

Now, applying our standard  truncations technique, see the detailed
discussion and bounds in \cite{Zeifman2014i,Zeifman2016t}, we can
find the limiting characteristics of the queue-length process, the
respective plots are shown in pictures \textcolor{red}{X-Y}.

\smallskip

{\bf Example 2.} Let $\mu_1=\mu_2=2,~~\la=3(1+\sin2\pi t),
~~\la_*=3$ (Example 3 from \cite{aa2018}). Put
$\delta=\sqrt{\frac{\mu}{\la_*}}=\frac{2}{\sqrt3}$. Then we obtain
$\alpha_1=2-3\epsilon(1+\sin2\pi t)$, $\alpha_3=5+3\sin2\pi
t-2\sqrt{3(1+\sin2\pi t)}$, $\alpha_4=7+3\sin2\pi
t+(4-\epsilon)\sqrt{3(1+\sin2\pi t)}.$

On the other hand, using Remark 2 we have simple corresponding
bounds: $\alpha_{10}=2-3\epsilon$, $\alpha_{30}=5 -2\sqrt{3}$,
$\alpha_{40}=7-(4-\epsilon)\sqrt{3}\ge7-4\sqrt{3}\ge0.07$. Hence
$\beta_{*0} =0.07.$

Thus we can obtain the following {\bf bound}
\begin{equation*}
\|\vp^*(t)-\vp^{**}(t)\|_1  \le
4\textsf{N}e^{-0.07t}\|\vz^*(0)-\vz^{**}(0)\|_{1\textsf{D}}.
\end{equation*}

\section{The case $\mu_1 >\mu_2$.}

First, let the intensities be constant. Put
$\mu_1(t)=(1+\chi)\mu_2(t)$, where  for $\chi >0.$ Then we have
$$\alpha_1=(1+\chi)\mu_2 -\epsilon \la,$$
$$\alpha_2=\la+\mu_2\left(1-\frac{\chi}{\epsilon}\right),$$
$$\alpha_3=\la(1-\delta_1)+(1+\chi)\mu_2,$$
$$\alpha_4=\la(1-\delta)+\mu_2\left(2+\chi-\frac{2+\epsilon+\chi}{\delta_1}\right),$$
and
$$\alpha_k=\la(1-\delta)+\mu_2\left(1-\frac{1}{\delta}\right)(2+\chi), \quad k \ge 5.$$

Put $\beta_*=\min_{i \le 4}(\alpha_i).$ Then we have
\begin{equation}
\gamma(B(t))_{1\textsf{D}}=-\min(\alpha_i(t))=-\beta_*.
\end{equation}

Hence, the following bound  on the rate of convergence holds:
\begin{equation}
\|\vp^*(t)-\vp^{**}(t)\|_1  \le 4e^{-\beta_*
t}\|\vz^*(0)-\vz^{**}(0)\|_{1\textsf{D}}.
\end{equation}

Let now the intensities $\la(t)$, $\mu_1(t) = (1+\chi)\mu_2(t)$ be
1-periodic. Put
$$\mu_{2*}=\int_0^1\mu_2(t)\,dt,~~\la_*=\int_0^1\la(t)\,dt.$$
Then, in accordance with Remark 2, we can find the corresponding
parameter $\beta_{0*}$ for the respective homogeneous model. Namely,
we have $\beta_{*0}=\min_{i\le 4} \left(\alpha_{i0}\right),$ and the
bound on the rate of convergence (\ref{beta*01}) for some positive
$\textsf{N}$.

\smallskip

{\bf Example 3.} Let $\mu_1(t)=6\left(1+\cos2\pi t\right)$,
$\mu_2(t)=5\left(1+\cos2\pi t\right)$, $\la(t) = 8\left(1+\sin2\pi
\right)$. (Example 1 from \cite{aa2018}). Then $\chi=0.2$,
$\mu_{2*}=5$, $\mu_{1*}=6$, $\la_*=8.$ Hence, we have
$$\alpha_{10}=6 -8\epsilon, \quad \alpha_{20}=13-\frac{1}{\epsilon},$$
$$\alpha_{30}=14 -8\delta_1, \quad  \alpha_{40}=19- 8\delta -\frac{11+5\epsilon}{\delta_1},$$
and
$$\alpha_{k0}=8(1-\delta)+11\left(1-\frac{1}{\delta}\right), \quad k \ge 5.$$

As we have already noted, the best value of the bound is attained,
when $\delta =\sqrt{\frac{\mu}{\la}}=\sqrt{\frac{11}{8}}$, then
$$\alpha_{k0}=\left(\sqrt{11} - \sqrt{8}\right)^2 \approx 0.2, \quad k \ge 5.$$

Now put $\epsilon=\frac{1}{12}$ and  $\delta_1=\frac{13}{8}$. Then
$\alpha_{10}>1$, $\alpha_{20}=1$,  $\alpha_{30}=1$ and
$\alpha_{40}=19- \sqrt{88} -\frac{11+5/12}{13/8} > 1$.

Then we obtain  $$\beta_{*0}= \inf_{i \ge 1}\alpha_{i0} =
\alpha_{50} = \left(\sqrt{11} - \sqrt{8}\right)^2 > 0.2,$$ and the
following  {\bf bound} on the rate of convergence holds:
\begin{equation*}
\|\vp^*(t)-\vp^{**}(t)\|_1 \le
4\textsf{N}e^{-0.2t}\|\vz^*(0)-\vz^{**}(0)\|_{1\textsf{D}}.
\end{equation*}

\subsubsection*{Acknowledgments.} The bounds on the rate of
convergence have been obtained by Zeifman and Korolev. The work of
Zeifman and Korolev is supported by the Russian Science Foundation
under grant 18-11-00155. Examples have been studied by Satin,
Kiseleva, Panfilova. 

\renewcommand{\refname}{References}

\clearpage

\begin{figure}
\begin{center}
  \includegraphics[width=14cm]{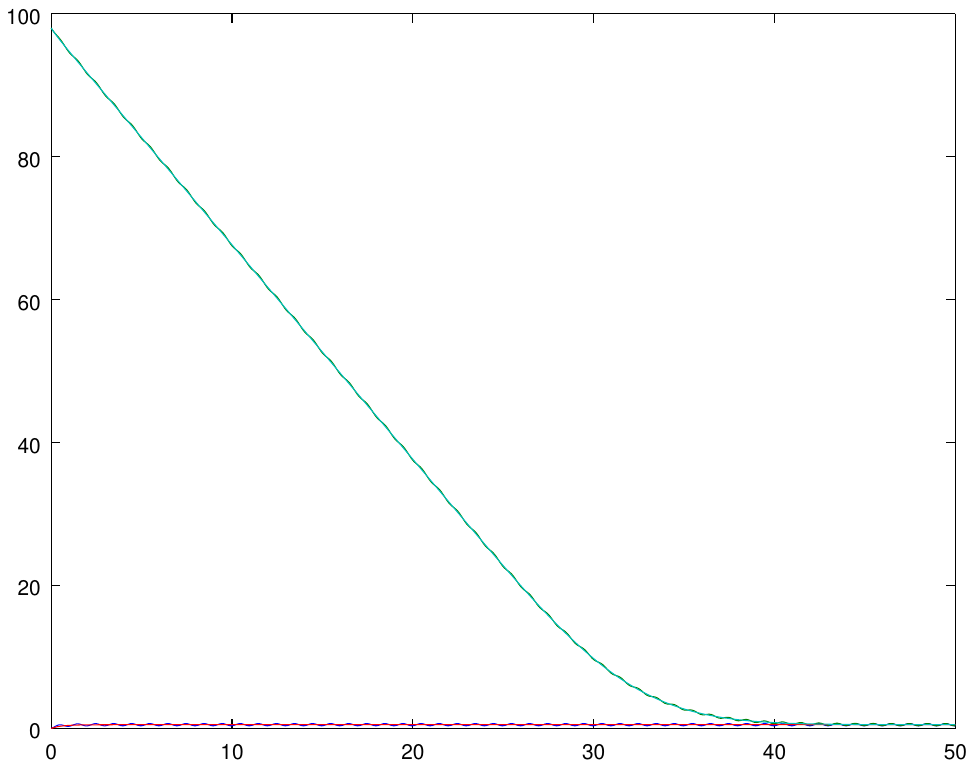}
\end{center}
\caption{Example 1. Approximation of the mean $E(t,k)$ for
$t\in[0,50]$ with initial conditions $X(0)=0$ and $X(0)=100$ for
original and homogeneous situations.}
\label{fig:08}       
\end{figure}

\begin{figure}
\begin{center}
  \includegraphics[width=14cm]{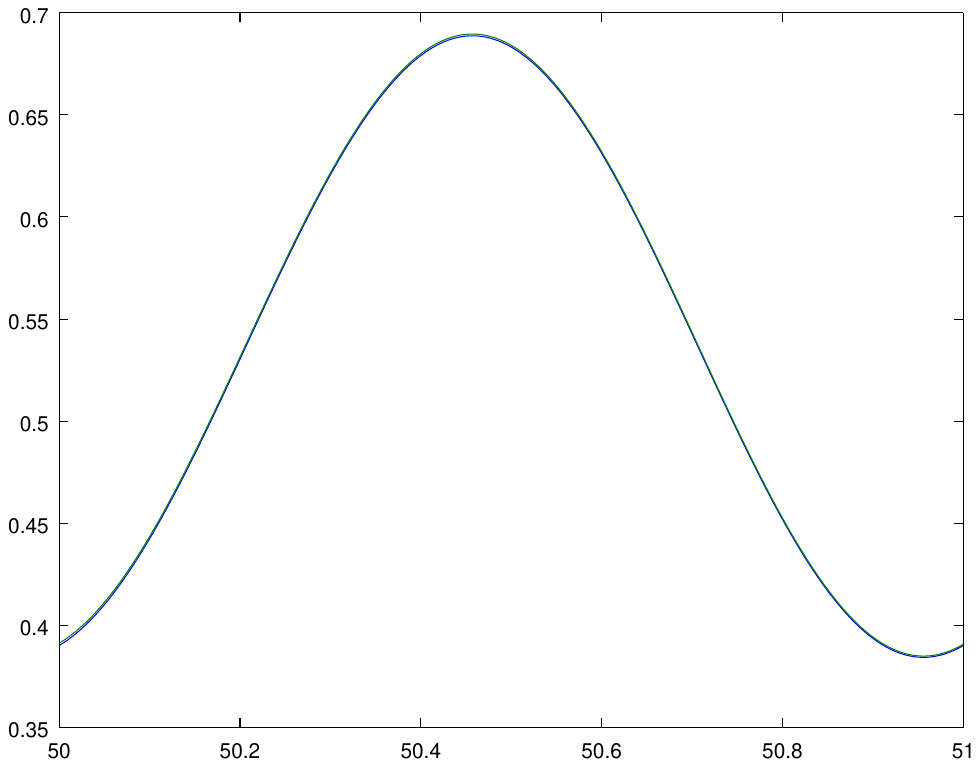}
\end{center}
\caption{Example 1. Approximation of the mean $E(t,k)$ for
$t\in[50,51]$ with initial conditions $X(0)=0$ and $X(0)=100$.}
\label{fig:09}       
\end{figure}

\begin{figure}
\begin{center}
  \includegraphics[width=14cm]{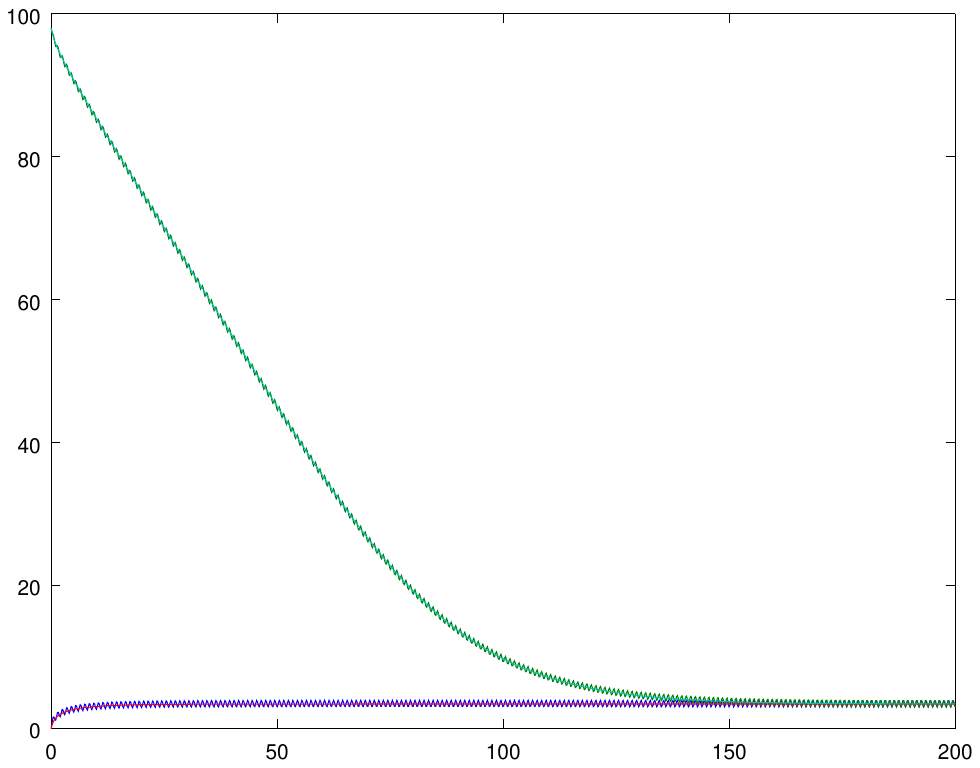}
\end{center}
\caption{Example 2. Approximation of the mean $E(t,k)$ for
$t\in[0,200]$ with initial conditions $X(0)=0$ and $X(0)=100$ for
original and homogeneous situations.}
\label{fig:18}       
\end{figure}

\begin{figure}
\begin{center}
  \includegraphics[width=14cm]{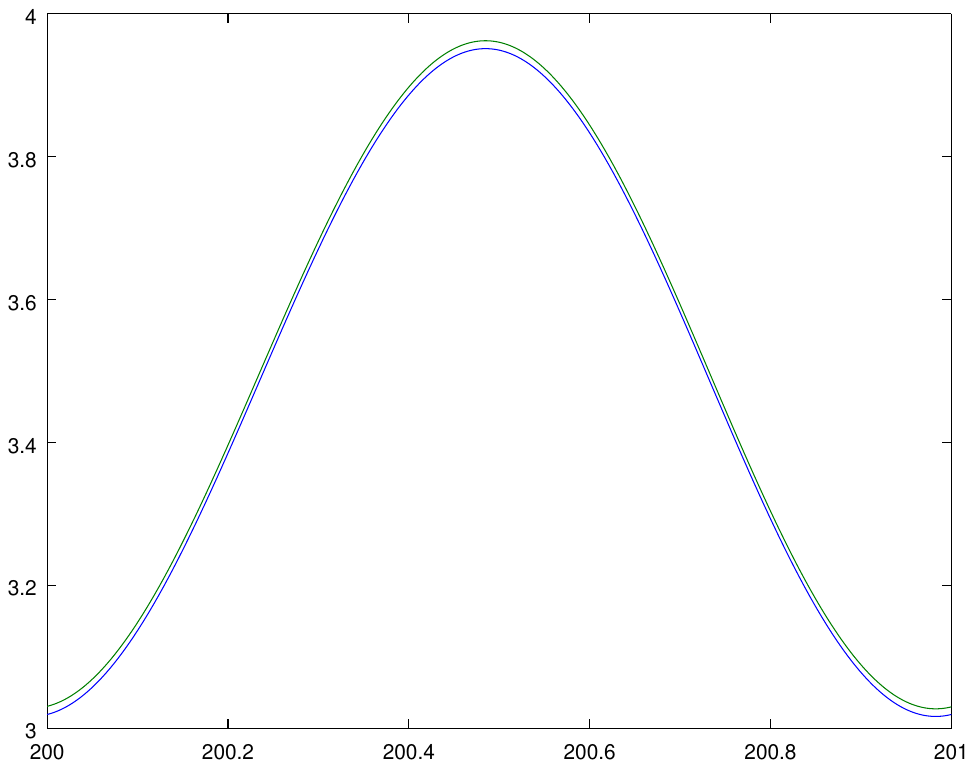}
\end{center}
\caption{Example 2. Approximation of the mean $E(t,k)$ for
$t\in[200,201]$ with initial conditions $X(0)=0$ and $X(0)=100$.}
\label{fig:19}       
\end{figure}

\begin{figure}
\begin{center}
  \includegraphics[width=14cm]{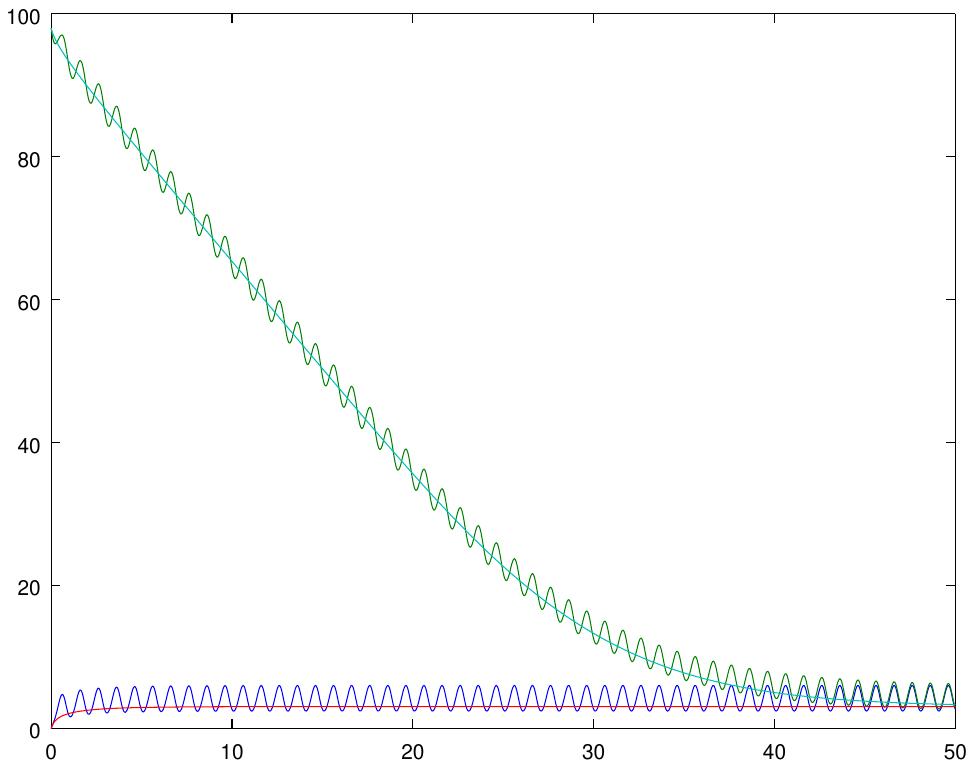}
\end{center}
\caption{Example 3. Approximation of the mean $E(t,k)$ for
$t\in[0,50]$ with initial conditions $X(0)=0$ and $X(0)=100$ for
original and homogeneous situations.}
\label{fig:38}       
\end{figure}

\begin{figure}
\begin{center}
  \includegraphics[width=14cm]{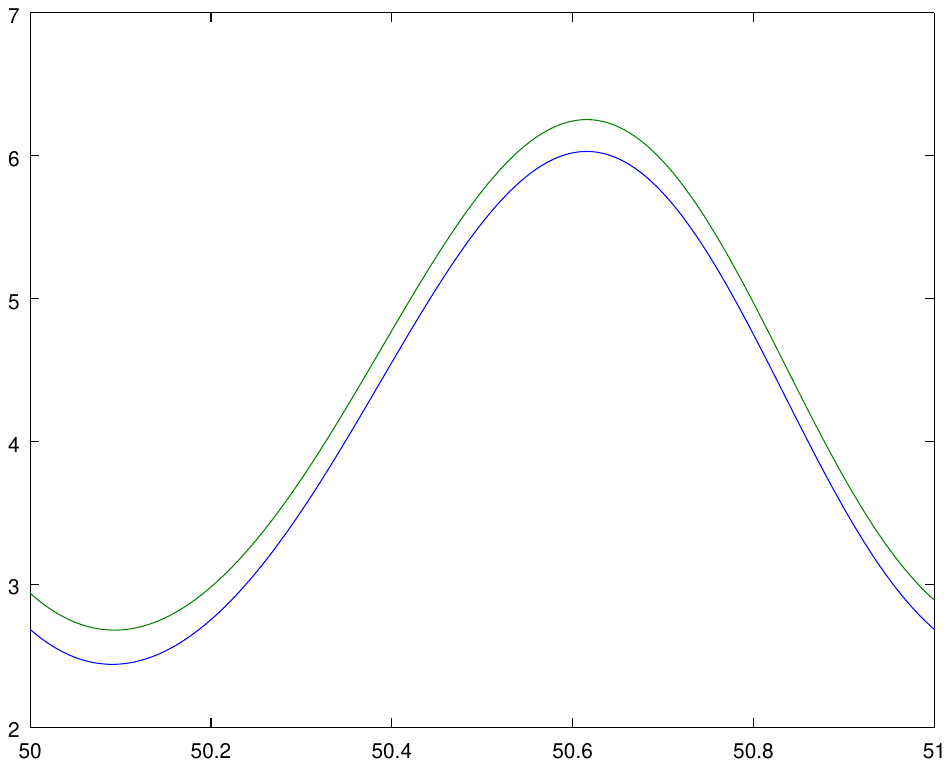}
\end{center}
\caption{Example 3. Approximation of the mean $E(t,k)$ for
$t\in[50,51]$ with initial conditions $X(0)=0$ and $X(0)=100$.}
\label{fig:39}       
\end{figure}

\smallskip

\begin{figure}
\begin{center}
  \includegraphics[width=14cm]{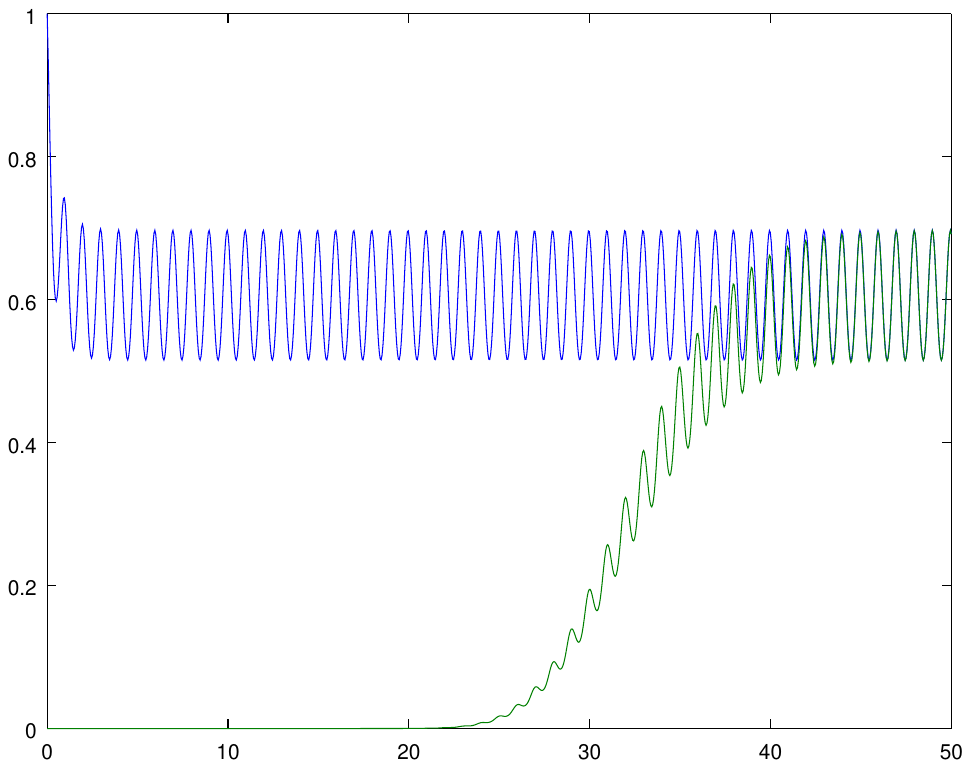}
\end{center}
\caption{Example 1. Approximation of the probability $P_{00}(t)$ for
$t\in[0,50]$ with initial conditions $X(0)=0$ and $X(0)=100$.}
\label{fig:01}       
\end{figure}


\begin{figure}
\begin{center}
  \includegraphics[width=14cm]{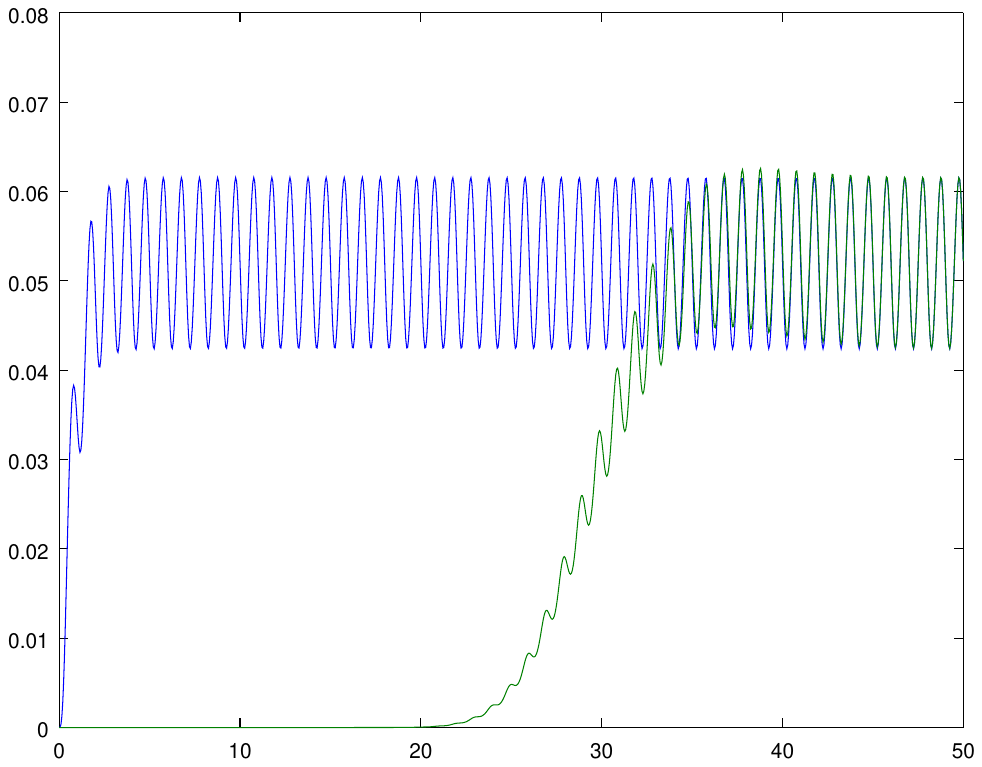}
\end{center}
\caption{Example 1. Approximation of the probability  $P_{01}(t)$
for $t\in[0,50]$ with initial conditions $X(0)=0$ and $X(0)=100$.}
\label{fig:02}       
\end{figure}


\begin{figure}
\begin{center}
  \includegraphics[width=14cm]{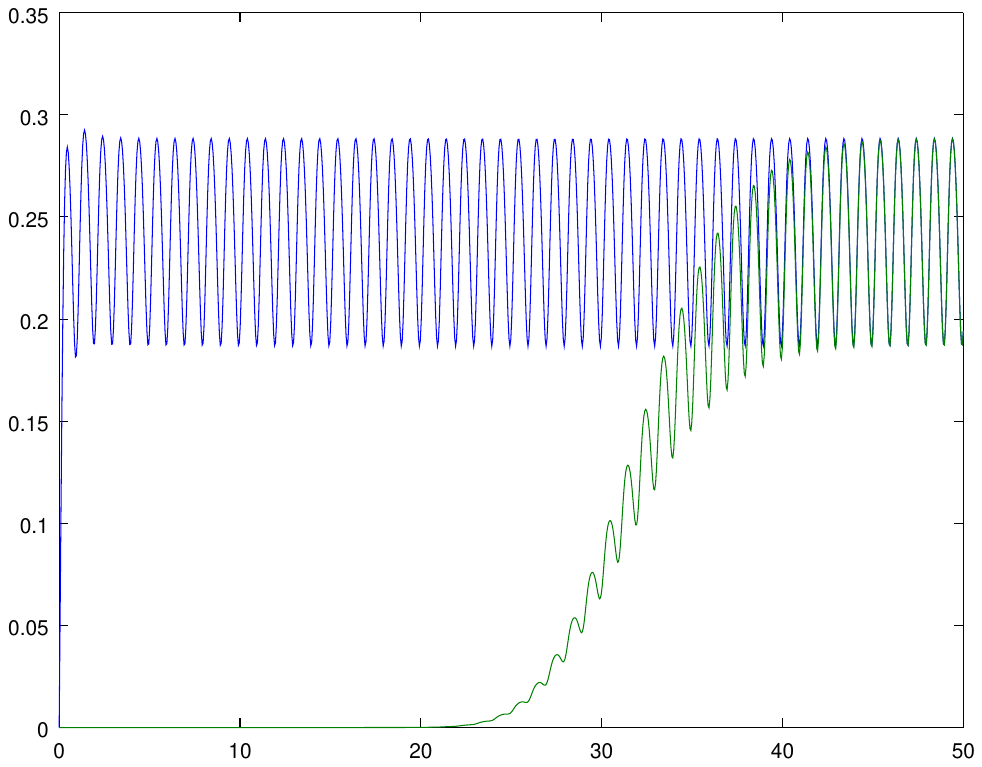}
\end{center}
\caption{Example 1. Approximation of the probability  $P_{10}(t)$
for $t\in[0,50]$ with initial conditions $X(0)=0$ and $X(0)=100$.}
\label{fig:03}       
\end{figure}


\begin{figure}
\begin{center}
  \includegraphics[width=14cm]{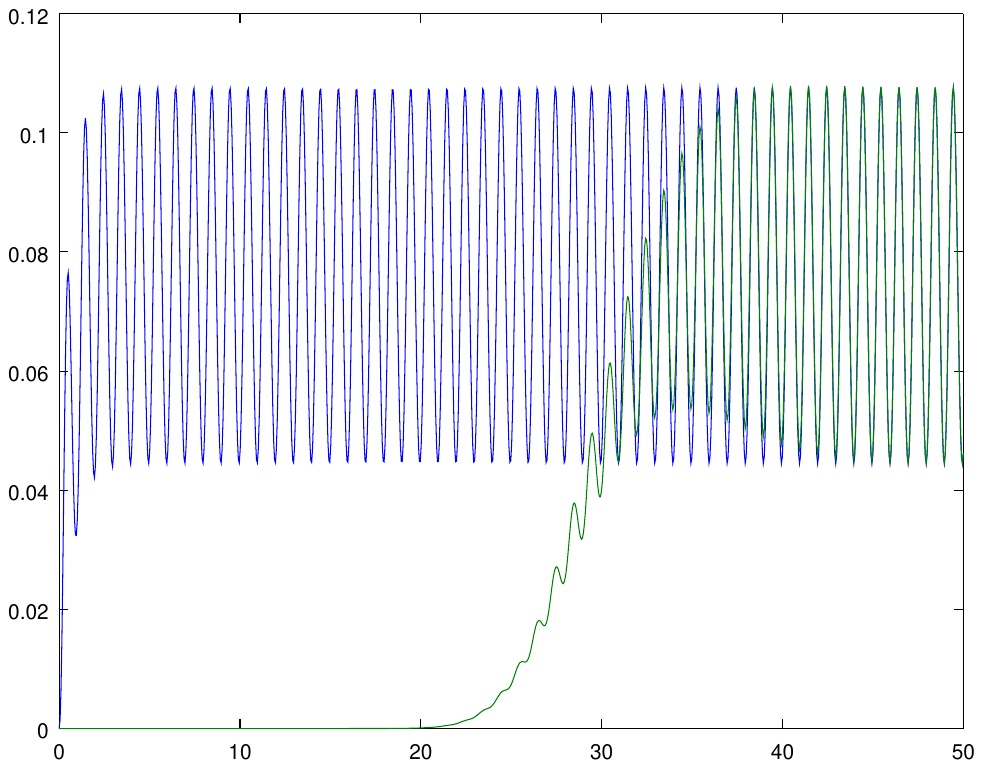}
\end{center}
\caption{Example 1. Approximation of the probability  $P_{11}(t)$
for $t\in[0,50]$ with initial conditions $X(0)=0$ and $X(0)=100$.}
\label{fig:04}       
\end{figure}


\begin{figure}
\begin{center}
  \includegraphics[width=14cm]{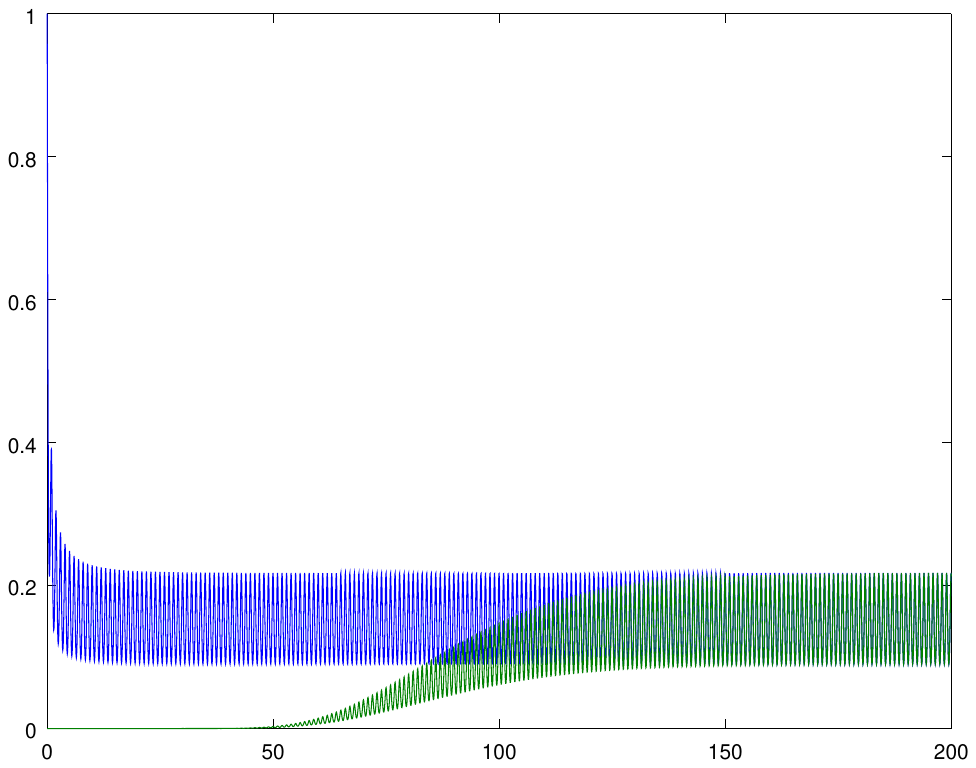}
\end{center}
\caption{Example 2. Approximation of the probability $P_{00}(t)$ for
$t\in[0,200]$ with initial conditions $X(0)=0$ and $X(0)=100$.}
\label{fig:11}       
\end{figure}

\begin{figure}
\begin{center}
  \includegraphics[width=14cm]{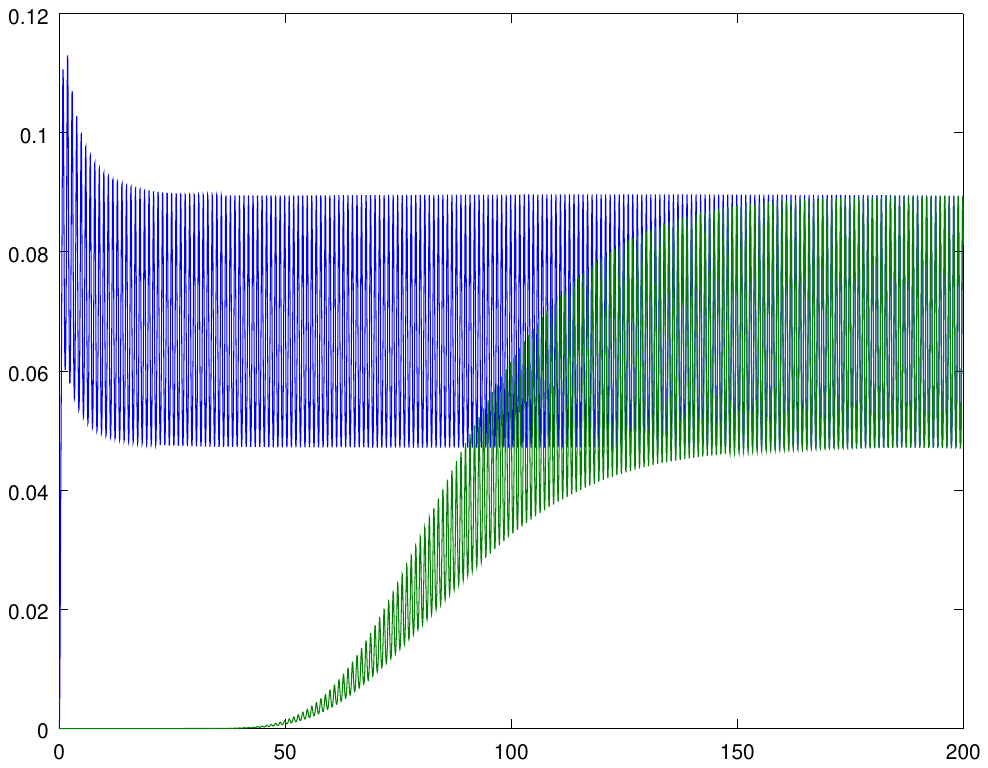}
\end{center}
\caption{Example 2. Approximation of the probability  $P_{01}(t)$
for $t\in[0,200]$ with initial conditions $X(0)=0$ and
$X(0)=100$.}
\label{fig:121}       
\end{figure}

\begin{figure}
\begin{center}
  \includegraphics[width=14cm]{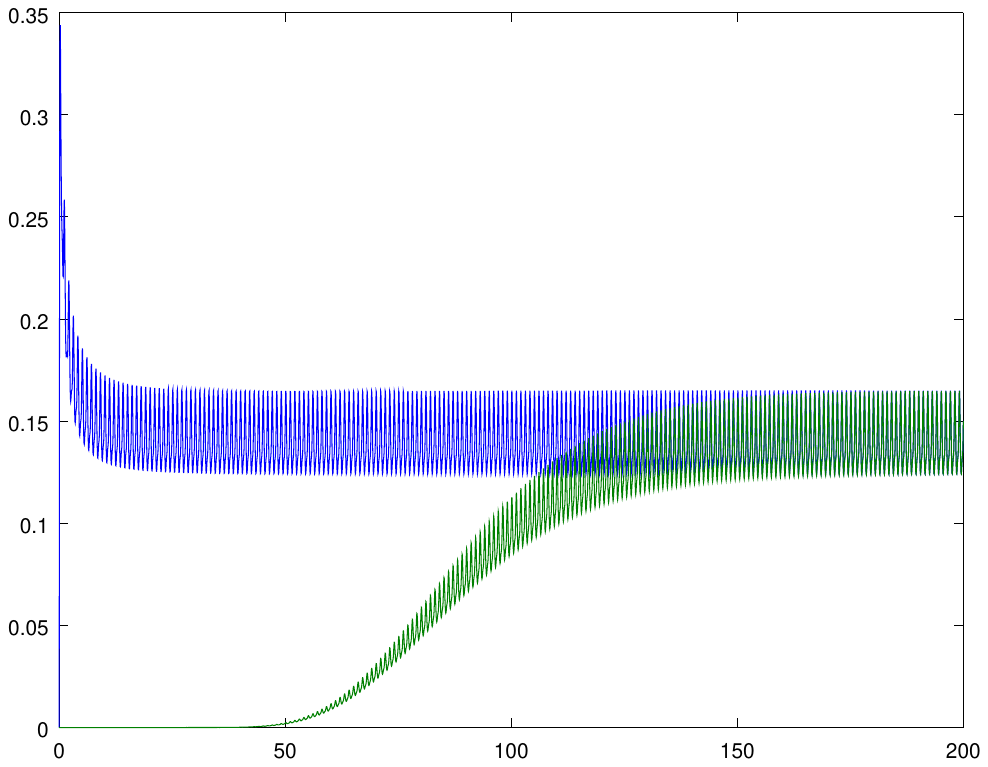}
\end{center}
\caption{Example 2. Approximation of the probability  $P_{10}(t)$
for $t\in[0,200]$ with initial conditions $X(0)=0$ and
$X(0)=100$.}
\label{fig:131}       
\end{figure}

\begin{figure}
\begin{center}
  \includegraphics[width=14cm]{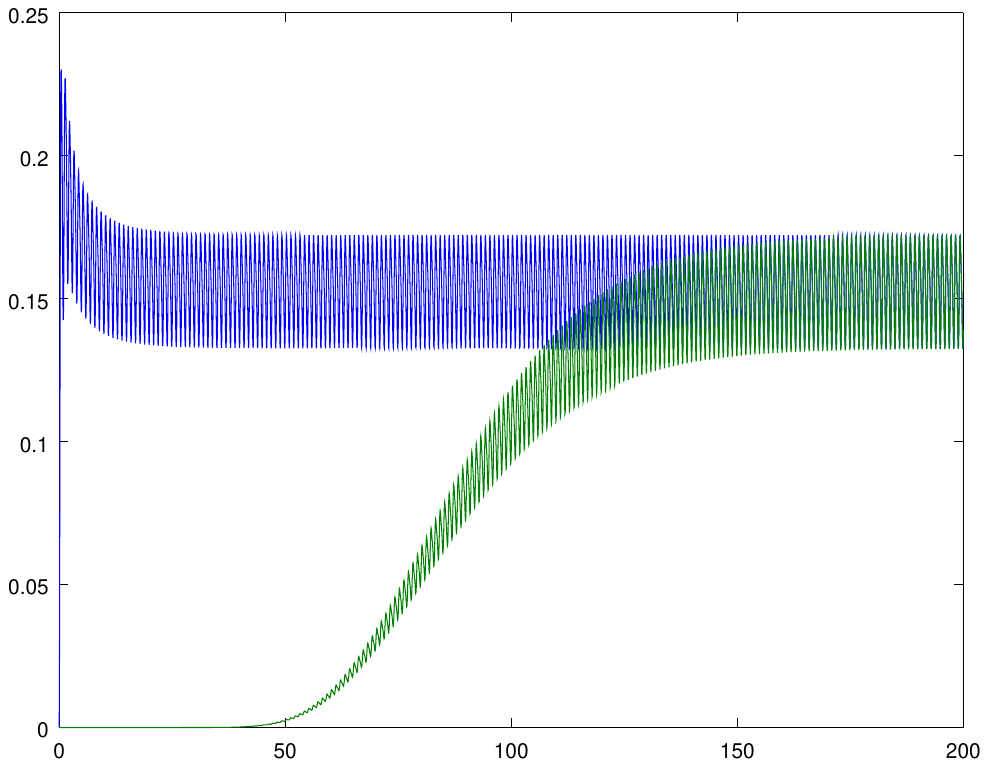}
\end{center}
\caption{Example 2. Approximation of the probability  $P_{11}(t)$
for $t\in[0,200]$ with initial conditions $X(0)=0$ and $X(0)=100$.}
\label{fig:141}       
\end{figure}

\begin{figure}
\begin{center}
  \includegraphics[width=14cm]{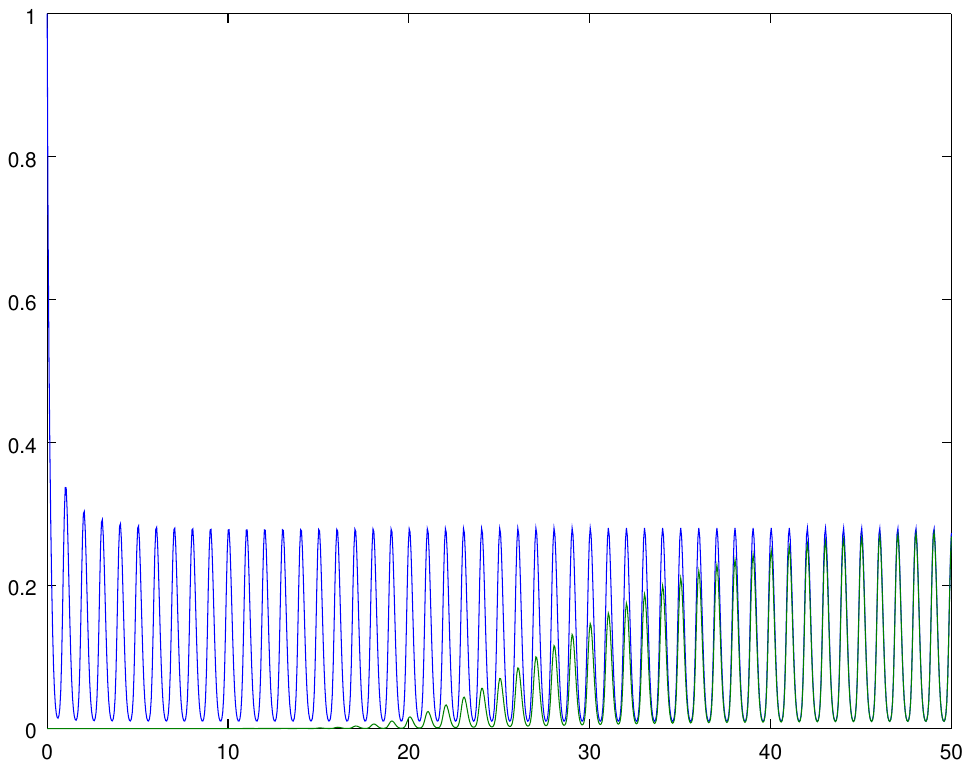}
\end{center}
\caption{Example 3. Approximation of the probability $P_{00}(t)$ for
$t\in[0,50]$ with initial conditions $X(0)=0$ and $X(0)=100$.}
\label{fig:31}       
\end{figure}

\begin{figure}
\begin{center}
  \includegraphics[width=14cm]{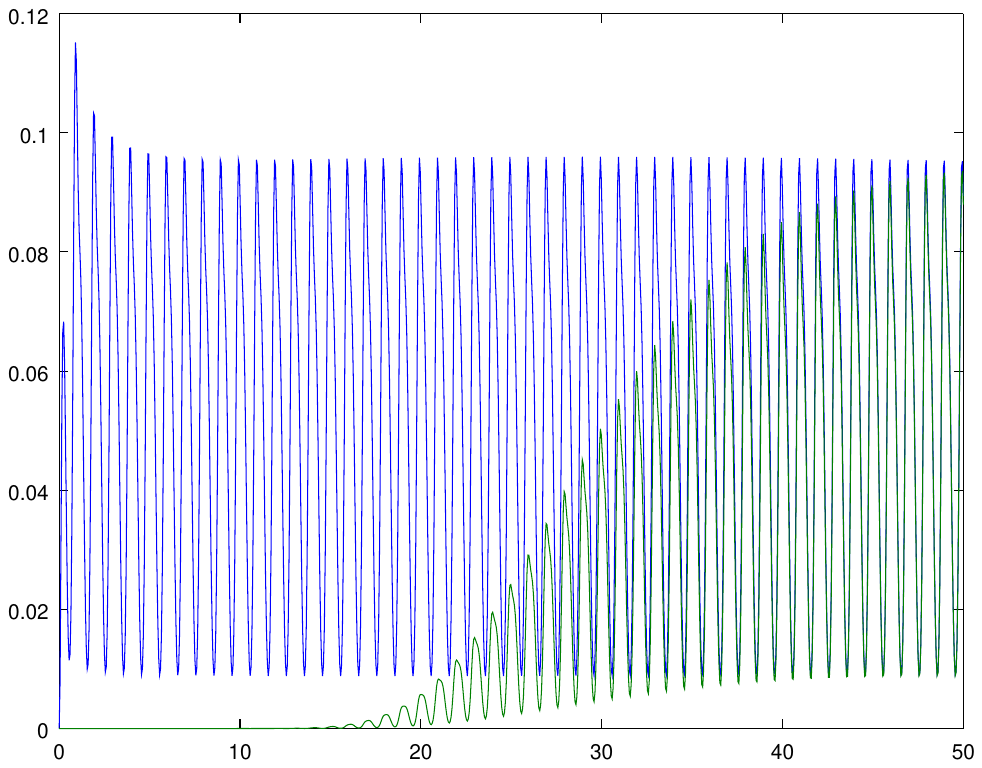}
\end{center}
\caption{Example 3. Approximation of the probability  $P_{01}(t)$
for $t\in[0,50]$ with initial conditions $X(0)=0$ and $X(0)=100$.}
\label{fig:32}       
\end{figure}

\begin{figure}
\begin{center}
  \includegraphics[width=14cm]{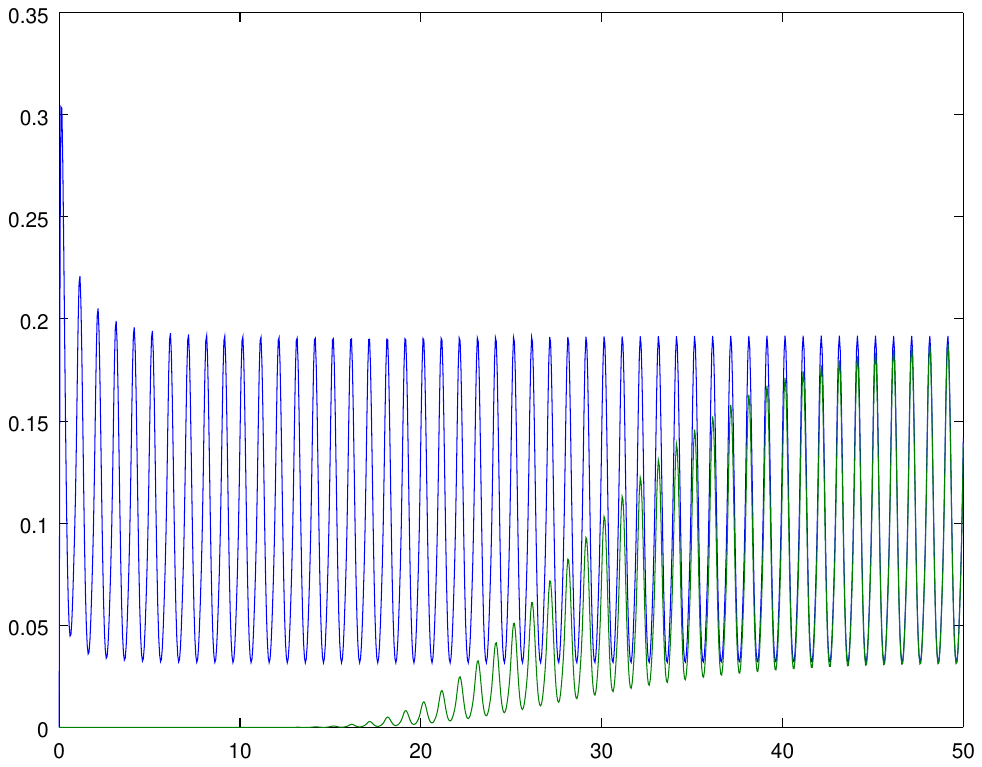}
\end{center}
\caption{Example 3. Approximation of the probability  $P_{10}(t)$
for $t\in[0,50]$ with initial conditions $X(0)=0$ and $X(0)=100$.}
\label{fig:33}       
\end{figure}

\begin{figure}
\begin{center}
  \includegraphics[width=14cm]{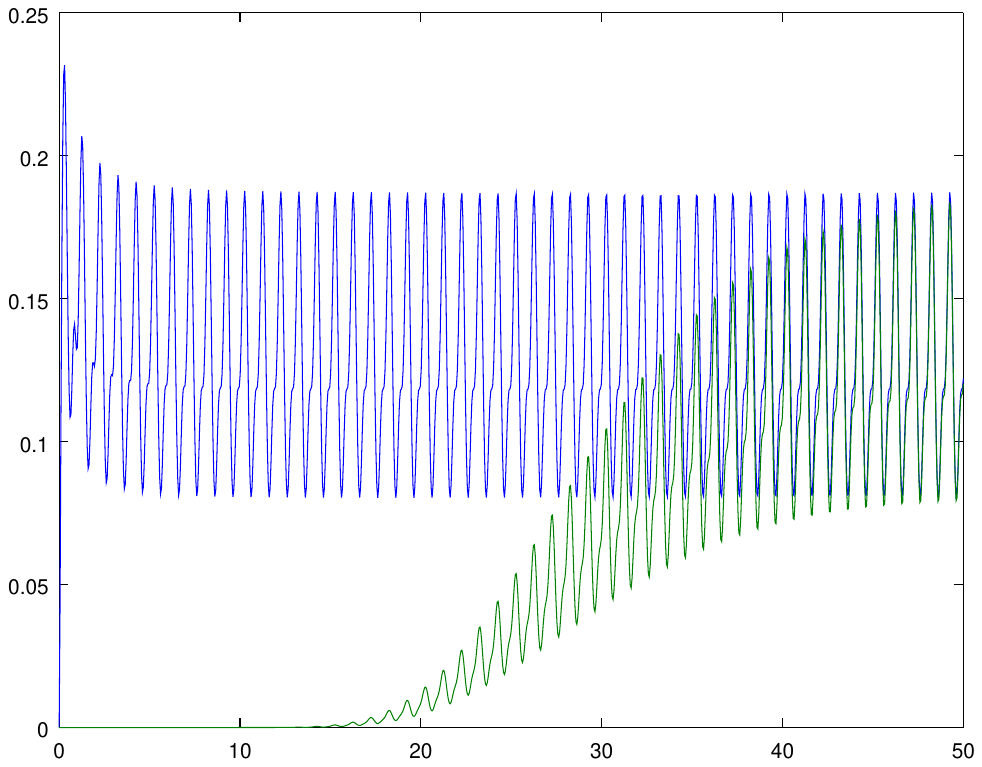}
\end{center}
\caption{Example 3. Approximation of the probability  $P_{11}(t)$
for $t\in[0,50]$ with initial conditions $X(0)=0$ and $X(0)=100$.}
\label{fig:34}       
\end{figure}

\end{document}